\documentclass[12pt]{amsart}%
\usepackage{amsmath}
\usepackage{graphicx}
\usepackage{amscd}
\usepackage{geometry}
\usepackage{amsfonts}
\usepackage{amssymb}%
\usepackage{amsthm,enumerate,graphicx,etex,diagrams,ifsym}
\newtheorem{theorem}{Theorem}[section]
\theoremstyle{plain}
\newtheorem{Theorem}{Theorem}[section]	
\theoremstyle{plain}

\newtheorem{Corollary}[theorem]{Corollary}

\newtheorem{Definition}{Definition}

\newtheorem{Lemma}[theorem]{Lemma}

\newtheorem{Problem}{Problem}

\newcommand{\BN}{\mathbb{N}}
\newcommand{\BZ}{\mathbb{Z}}

\newcommand{\BS}{\mathbb{S}}

\newenvironment{Proof}{\par\noindent{\sc Proof}\quad}{\hfill\qed\par\smallskip}
\newenvironment{RestateTheorem}[3]{\par\vspace{12pt}\noindent{\bf #1~\ref{#2}}(#3){\bf .}\it}{\par\vspace{12pt}}
\let\oldtocsection=\tocsection
\let\oldtocsubsection=\tocsubsection
\let\oldtocsubsubsection=\tocsubsubsection
\renewcommand{\tocsection}[2]{\hspace{0em}\oldtocsection{#1}{#2}}
\renewcommand{\tocsubsection}[2]{\hspace{2em}\oldtocsubsection{#1}{#2}}
\renewcommand{\tocsubsubsection}[2]{\hspace{4em}\oldtocsubsubsection{#1}{#2}}

\begin{document}

\title[More Examples of Pseudo-Collars]{More Examples of Pseudo-Collars on High-Dimensional Manifolds}
\author{Jeffrey J. Rolland }
\address{Department of Mathematics, Milwaukee School of Engineering,
Milwaukee, Wisconsin 53202}
\email{rolland@msoe.edu}
\date{March 22, 2018}

\begin{abstract}

In a previous paper, we developed general techniques for constructing a variety of pseudo-collars, as defined by Guilbault and Tinsley, with roots in earlier work by Chapman and Siebenmann. As an application of our techniques, we exhibited an uncountable collection of pseudo-collars, all with the same boundary and similar fundamental group systems at infinity. Construction of that family was very specific; it relied on properties of Thompson's group $V$. In this paper, we provide a more general approach to constructing similar collections of examples. Instead of using Thompson's group $V$, we base our new examples on a broader and more common collection of groups, in particular, fundamental groups of certain hyperbolic 3-manifolds.

\end{abstract}
\maketitle
%\tableofcontents

\section{Introduction and Main Result \label{Section: Introduction and Main Result}}

%Several lines of investigation in geometric topology. The recognition
%criterion follows.

In \cite{Rolland1}, a geometric procedure for producing a ``reverse'' to Quillen's plus construction, a construction called a \textit{1-sided h-cobordism} or \textit{semi-h-cobordism}, was developed. This reverse to the plus construction was then used to produce uncountably many distinct ends of manifolds called \textit{pseudo-collars}, which are stackings of 1-sided h-cobordisms. The notion of pseudo-collars originated in \cite{C-S} in Hilbert cube manifold theory, where it was part of a necessary and sufficient condition for placing a $\mathcal{Z}$-set as the boundary of an open Hilbert cube manifold. In \cite{Guilbault1}, \cite{G-T2}, and \cite{G-T3}, the authors were interested in pseudo-collars on finite-dimensional manifolds for the same reason, attempting to put a $\mathcal{Z}$-set as the boundary of an open high-dimensional manifold. Each of the pseudo-collars had the same boundary and pro-homology systems at infinity and similar group-theoretic properties for their pro-fundamental group systems at infinity. In particular, the kernel group of each group extension for each 1-sided h-cobordism in the pseudo-collars was the same group. Nevertheless, the pro-fundamental group systems at infinity were all distinct. A good deal of combinatorial group theory was needed to verify this fact, including an application of Thompson's group V. 

\indexspace

In this paper, we extend this construction to have the kernel groups of each group extension be a free product $S*S$ of any finitely presented, centerless, superperfect group $S$ which contains a countably infinite list of elements $\{a_1, a_2, a_3, \ldots\}$ with the property that for any isomorphism $\phi: S \rightarrow S$, $\phi(a_i) \ne a_j$ for $i \ne j$. Note that this class of groups includes the fundamental group of any closed hyperbolic manifold $M$ with first two homology groups $H_1(M) = H_2(M) = 0$. Note further that the $(1,n)$ Dehn filling of the figure-8 knot complement in $\mathbb{S}^3$ for $n > 1$ satisfies this condition, so already we have a countably infinite collection of such kernel groups, for each of which we produce an uncountable collection of distinct pseudo-collars.

\indexspace

We work in the category of smooth manifolds, but all our results apply equally well to the categories of PL and topological manifolds. The manifold version of Quillen's plus construction provides a way of taking a closed, smooth manifold $M$ of dimension $n \ge 5$ whose fundamental group $G = \pi_1(M)$ contains a perfect normal subgroup $P$ which is the normal closure of a finite number of elements and producing a compact cobordism $(W,M,M^+)$ to a manifold $M^+$ whose fundamental group is isomorphic to $Q = G/P$ and for which $M^+ \hookrightarrow W$ is a simple homotopy equivalence. By duality, the map $f:M \rightarrow M^+$ given by including $M$ into $W$ and then retracting onto $M^+$ induces an isomorphism $f_*:H_*(M;\mathbb{Z}Q) \rightarrow H_*(M^+;\mathbb{Z}Q)$ of homology with twisted coefficients. By a clever application of the s-Cobordism Theorem, such a cobordism is uniquely determined by $M$ and $P$ (see \cite{Freedman-Quinn} P. 197).

\indexspace

In ``Manifolds with Non-stable Fundamental Group at Infinity I'' \cite{Guilbault1}, Guilbault outlines a structure to put on the ends of an open smooth manifold $N$ with finitely many ends called a \textit{pseudo-collar}, which generalizes the notion of a collar on the end of a manifold introduced in Siebenmann's dissertation \cite{Siebenmann}. A pseudo-collar is defined as follows. Recall that a manifold $U^n$ with compact boundary is an open collar if $U^n \approx \partial U^n \times [0,\infty)$; it is a homotopy collar if the inclusion $\partial U^n \hookrightarrow U^n$ is a homotopy equivalence. If $U^n$ is a homotopy collar which contains arbitrarily small homotopy collar neighborhoods of infinity, then we call $U^n$ a \textit{pseudo-collar}. We say that an open $n$-manifold $N^n$ is collarable if it contains an open collar neighborhood of infinity, and that $N^n$ is \textit{pseudo-collarable} if it contains a pseudo-collar neighborhood of infinity.

\indexspace

Each pseudo-collar admits a natural decomposition as a sequence of compact cobordisms $(W,M,M_-)$, where $W$ is a 1-sided h-cobordism (see Definition \ref{defsemi-h-cob} below). If a 1-sided h-cobordism is actually an s-cobordism (again, see Definition \ref{defsemi-h-cob} below), it follows that the cobordism $(W,M_-,M)$ is a a plus cobordism. (This somewhat justifies the use of the symbol ``$M_-$'' for the right-hand boundary of a 1-sided h-cobordism, a play on the traditional use of $M^+$ for the right-hand boundary of a plus cobordism.)

\indexspace

The general problem of a reverse to Quillen's plus construction in the high- \\ dimensional manifold category is as follows.

\begin{Problem}[Reverse Plus Problem]
Suppose $G$ and $Q$ are finitely-presented groups and $\Phi: G \twoheadrightarrow Q$ is an onto homomorphism with $\ker(\Phi)$ perfect. Let $M^n$ ($n \ge 5$) be a closed smooth manifold with $\pi_1(M) \cong Q$.

\indexspace

Does there exist a compact cobordism $(W^{n+1}, M, M_-)$ with 

\begin{diagram}[size=14.5pt]
1 & \rTo & \ker(\iota_{\#}) & \rTo & \pi_1(M_-) & \rTo^{\iota_{\#}} & \pi_1(W) & \rTo & 1
\end{diagram}

equivalent to

\begin{diagram}[size=14.5pt]
1 & \rTo & \ker(\Phi) & \rTo & G & \rTo^{\Phi} & Q & \rTo & 1
\end{diagram}

and $M \hookrightarrow W$ a (simple) homotopy equivalence. 
\end{Problem}

Notes:
\begin{itemize}

\item The fact that $G$ and $Q$ are finitely presented forces $\ker(\Phi)$ to be the normal closure of a finite number of elements. (See, for instance, \cite{Guilbault1} or \cite{Siebenmann}.)

\item Closed manifolds $M^n$ ($n \ge 5$) in the various categories (smooth, PL, and topological) with $\pi_1(M)$ isomorphic to a given finitely presented group $Q$ always exist. In the smooth category, one starts with an $n$-disk $\mathbb{B}^n$, attached one 1-handle for each generator of a given presentation of $Q$, and then attaches one 2-handle for each relator of the given presentation of $Q$, following a path in the 1-handles of $\mathbb{S}^n$ with the 1-handles attached for the corresponding relator for the attaching map of the 2-handle. The boundary is a closed manifold $M$ with desired fundamental group $Q$. Similar procedures exist in the other categories.

\end{itemize}

\begin{Definition}\label{defsemi-h-cob}
Let $N^n$ be a compact smooth manifold. A \textbf{1-sided h-cobordism} $(W,N,M)$ is a cobordism with either $N \hookrightarrow W$ or $M \hookrightarrow W$ is a homotopy equivalence (if it is a simple homotopy equivalence, we call $(W,N,M)$ a \textbf{1-sided s-cobordism}). [A 1-sided h-cobordism  $(W,N,M)$ is so-named presumably because it is ``one side of an h-cobordism''].
\end{Definition}

One wants to know under what circumstances 1-sided h-cobordisms exists, and, if they exist, how many there are. Also, one is interested in controlling the torsion and seeing when it can be eliminated.

\indexspace

There are some cases in which 1-sided h-cobordisms are known not to exist. For instance, if $P$ is finitely presented and perfect but not superperfect, $Q = \langle e \rangle$, and $M = \BS^n$, then a solution to the Reverse Plus Problem would produce an $M_-$ that is a homology sphere. But it is a standard fact that a manifold homology sphere must have a a superperfect fundamental group! (See, for instance, \cite{Kervaire}.) (A perfect group $P$ may be described as a group with $H_1(P) = 0$. A \emph{superperfect} group $S$ is then one with $H_1(S) = H_2(S) = 0$.)

\indexspace

The key point is that the solvability of the Reverse Plus Problem depends not just upon the group data, but also upon the manifold $M$ with which one begins. For instance, one could start with a group $P$ which is finitely presented and perfect but not superperfect, let $N_-$ be a manifold obtained from the boundary of a regular neighborhood of the embedding of a presentation 2-complex for $P$ in $\BS^{n+1}$, and let $(W, N_-, N)$ be the result of applying Quillen's plus construction to to $N_-$ with respect to all of $P$. Then again $Q = \langle e \rangle$ and $\Phi: P \twoheadrightarrow Q$ but $N$ clearly admits a 1-sided s-cobordism, namely $(W, N, N_-)$ (however, of course, we cannot have $N$ a sphere or $N_-$ a homology sphere).

\indexspace

An underlying goal of papers \cite{Guilbault1}, \cite{G-T2}, and \cite{G-T3} is to understand when non-compact manifolds with compact (possibly empty) boundary admit $\mathcal{Z}$-compactifications. In \cite{C-S}, it is shown that a Hilbert cube manifold admits a $\mathcal{Z}$-compactification if and only if it is pseudo-collarable and the Whitehead torsion of the end can be controlled. In \cite{Guilbault2}, Guilbault asks whether the universal cover of a closed, aspherical manifold ($n \ge 6$) is always pseudo-collarable. He further asks if pseudo-collarability plus control of the Whitehead torsion of the end is enough for finite-dimensional manifolds ($n \ge 6$) to admit a $\mathcal{Z}$-compactification. Still further, he shows in \cite{Guilbault2} that any two $\mathcal{Z}$-boundaries of an ANR must be shape equivalent. Finally, he and Ancel show in \cite{A-G} that if two closed, contractible manifolds $M^n$ and $N^n$ ($n \ge 6$) admit homeomorphic boundaries, then $M$ is homeomorphic to $N$. This is most interesting when the contractible manifolds are universal covers of closed aspherical manifolds. In that case, these questions may be viewed as an approach to the famous Borel Conjecture, which asks whether two aspherical manifolds with isomorphic fundamental group are necessarily homeomorphic.

\indexspace

In this paper, we are interested in producing families of pseudo-collars on high-dimensional closed manifolds with fundamental group $\BZ$, but we aim to use more kernel groups than than are used in \cite{Rolland1} -- just the free product of two copies of Thompson's group $V$ as the kernel group -- for each group extension.

\indexspace

Here is a statement of our main result.

\begin{Theorem}[Uncountably Many Pseudo-Collars on Closed Manifolds with the Same Boundary and Similar Pro-$\pi_1$]  \label{thmpseudo-collars}
Let $M^n$ be a closed smooth manifold ($n \ge 6$) with $\pi_1(M) \cong \BZ$ and let $S$ be the finitely presented group $K*K$ which is the free product of 2 copies of an admissible group (defined below) group $K$. Then there exists an uncountable collection of pseudo-collars $\{N^{n+1}_{\omega}\ |\ \omega \in \Omega\}$, no two of which are homeomorphic at infinity, with  $\partial N^{n+1}_{\omega} = M^n$, each with fundamental group at infinity that may be represented by an inverse sequence

\begin{diagram}[size=14.5pt]
\BZ & \lOnto^{\alpha_1} & G_{1} & \lOnto^{\alpha_2} & G_{2} & \lOnto^{\alpha_3} & G_{3} & \lOnto^{\alpha_4} & \ldots \\
\end{diagram}

with $\ker(\alpha_i) = S$ for all $i$.
\end{Theorem}

The authors would like to thank Craig Guilbault, Derek Holt, Jason Manning, and Marston Conder for helpful conversations.

\section{Proof of the Main Result \label{Section: Proof of the Main Result}}
\begin{Definition}
We call a group $K$ \textbf{admissible} if it is a finitely presented, centerless, Hopfian, co-Hopfian, freely indecomposable, superperfect group which admits a countably infinite list of elements $A = \{a_1, a_2, a_3, \ldots\}$, called an \textbf{unpermutable subset}, with the property that for any isomorphism $\phi: K \rightarrow K$, $\phi(a_i) \ne a_j^{\pm 1}$ for $i \ne j$. 
\end{Definition}

Note that this class of groups includes the fundamental group of any closed, orientable hyperbolic manifold $M$ with first two homology groups $H_1(M) = H_2(M) = 0$. They are centerless as follows. Let $K$ be the fundamental group of a hyperbolic manifold. Let $ a \ne e \ne b$ be in $K$ with neither $a$ nor $b$ a power of some common element $c$. Then $a$ and $b$ don't commute. For, suppose momentarily that they did commute. Then, as the fundamental group of a hyperbolic manifold do not admit torsion (see, for example, Proposition 2.45 in \cite{Hatcher}), they would generate a $\BZ \oplus \BZ$ subgroup of $K$, which is impossible (see, for example, Theorem 7.1 The Flat Torus Theorem in \cite{B-H}). They are all Hopfian by a theorem of Sela in \cite{Sela}. They are co-Hopfian again by a theorem of Sela in \cite{Sela2}.  By Stallings's Theorem in \cite{Stallings}, they are all freely indecomposable as they all have one end. They all satisfy the condition on having an unpermutable subset as follows. If $a$ is any generator in a presentation of its fundamental group, any element of $\{a^k\ |\ k \in \BZ \backslash \{0\}\}$ has the length of a geodesic representative in the natural hyperbolic Riemanninan metric as an isometry invariant, so, by Mostow Rigitity, this is a countably infinite list of elements with the property that for any isomorphism $\phi: K \rightarrow K$, $\phi(a_i) \ne a_j^{\pm 1}$ for $i \ne j$. Note further that the fundamental group of the $(1,q)$ Dehn filling of the figure-8 knot complement in $\mathbb{S}^3$ for $q > 1$ is a superperfect group which is the fundamental group of a hyperbolic manifold, and therefore has all the required properties of the preceding paragraph. This is perfect as $p=1$ in the Dehn filling, so the Dehn filling kills $\pi_1(M)$. Since $M$ is an closed, orientable 3-manifold, it is superperfect by universal coefficients and Poincar\`e duality. Thus, already we have a countably infinite collection of admissible groups, for each of which we produce an uncountable collection of distinct pseudo-collars.

\indexspace

Let $K$ be an admissible group and $A$ its unpermutable subset, order $A$ as \\ $(a_1, a_2, a_3, \ldots)$, and let $a_i$ denote the $i^{th}$ element of $A$. Let $S = K_1*K_2$, where each $K_i$ is a copy of the above-mentioned admissible group.

\indexspace

Recall, if $H$ is a group, $Aut(H)$ is the automorphism group of $H$. Define $\mu: H \rightarrow Aut(H)$ to be $\mu(h)(h') = hh'h^{-1}$. Then the image of $\mu$ in $Aut(H)$ is called \textit{the inner automorphism group of $H$}, $Inn(H)$. The inner automorphism group of a group $H$ is always normal in $Aut(H)$. The quotient group $Aut(H)/Inn(H)$ is called the \textit{outer automorphism group} $Out(H)$. The kernel of $\mu$ is called the \textit{center of $H$}, $Z(H)$; it is the set of all $h \in H$ such that for all $h' \in H, hh'h^{-1} = h'$. One has the exact sequence

\indexspace

\begin{diagram}
1 & \rTo & Z(H) & \rTo  & H & \rTo^\mu  & Aut(H) & \rTo^\alpha & Out(H) & \rTo & 1 \\
\end{diagram}

\indexspace

Define a map $\Phi: K_1 \rightarrow Out(K_1*K_2)$ by $\Phi(a) = \phi_a$, where $\phi_a \in Out(K_1*K_2)$ is the outer automorphism defined by the automorphism 

\indexspace

$\phi_{a}(k) = 
\begin{cases}
k & \text{if } k \in K_1 \\
aka^{-1} & \text{if } k \in K_2 \\
\end{cases}$

\indexspace

($\phi_a$ is called a \emph{partial conjugation}.)

Lemma \ref{lemnontriv-free-prod} and Corollaries \ref{cornontriv-free-product}, \ref{lemstraightening-up-lemma}, and \ref{corstraightening-up-corollary} first appeared in \cite{Rolland1}; the reader is referred there for proofs.

\begin{Lemma} \label{lemnontriv-free-prod}
Let $A, B, C, \text{and } D$ be non-trivial groups. Let $\phi: A \times B \rightarrow C*D$ be a surjective homomorphism. Then one of $\phi(A \times \{1\})$ and $\phi(\{1\} \times B)$ is trivial and the other is all of $C*D$
\end{Lemma}

\begin{Corollary}  \label{cornontriv-free-product}
Let $A_1, \ldots, A_n$ be non-trivial groups and let $C*D$ be a free product of non-trivial groups. Let $\phi: A \times \ldots \times A_n \rightarrow C*D$ be a surjective homomorphism.

\indexspace

Then one of the $\phi(\{1\} \times \ldots A_i \times \ldots \times \{1\})$ is all of $C*D$ and the rest are all trivial.

\end{Corollary}

\begin{Corollary} \label{lemstraightening-up-lemma}
Let $S_1, S_2, \ldots, S_n$ all be copies of the same non-trivial free product, and let $\psi: S_1 \times S_2 \times \ldots \times S_n \rightarrow S_1 \times S_2 \times \ldots \times S_n$ be a isomorphism. Then $\psi$ decomposes as a ``matrix of maps'' $\psi_{i,j}$, where each $\psi_{i,j} = \pi_{S_j} \circ \psi|_{S_i}$ (where $\pi_{S_j}$ is projection onto $S_j$), and there is a permutation $\sigma$ on $n$ indices with the property that each $\psi_{\sigma(j), j}: S_{\sigma(j)} \rightarrow S_j$ is an isomorphism, and all other $\psi_{i,j}$'s are the zero map.
\end{Corollary}

Note that the $\psi_{i,j}$'s form a matrix where each row and each column contain exactly one isomorphism, and the rest of the maps are trivial maps - what would be a permutation matrix (see page 100 in \cite{Robbin}, for instance) if the isomorphisms were replaced by ``1'''s and the trivial maps were replaced by ``0'''s.

\begin{Corollary} \label{corstraightening-up-corollary}
Let $S_1, S_2, \ldots, S_n$ all be copies of the same non-trivial Hopfian free product, and let $\psi: S_1 \times S_2 \times \ldots \times S_n \rightarrow S_1 \times S_2 \times \ldots \times S_m$ be a epimorphism with $m < n$. Then $\psi$ decomposes as a ``matrix of maps'' $\psi_{i,j} = \pi_{S_j} \circ \psi|_{S_i}$, and there is a 1-1 function $\sigma$ from the set $\{1, \ldots, m\}$ to the set $\{1, \ldots, n\}$ with the property that $\psi_{\sigma(j), j}: S_{\sigma(j)} \rightarrow S_{j}$ is an isomorphism, and all other $\psi_{i,j}$'s are the zero map.
\end{Corollary}

Lemma \ref{lemconderisomorphism-lemma} and its corollary, Lemma \ref{lemconderepimorphism-lemma}, are the only lemmas that differ principally from the exposition of \cite{Rolland1} and constitute the main new material for this paper. We give their proofs in some detail.

\begin{Lemma} \label{lemconderisomorphism-lemma}
Let $K$ be an admissible group and $A$ one of its unpermutable subsets, and let $S = K_1*K_2$, where each $K_i$ is a copy of $K$. Let $S_1, \ldots S_n$ be $n$ copies of $S$, let $a_1, a_2, \ldots, a_n$ be distict elements of $A$, and define $\phi_{(a_1, \ldots, a_n)}: S_1 \times \ldots \times S_n \rightarrow S_1 \times \ldots S_n$ by $\phi_{(a_1, \ldots a_n)}(x_1, \ldots, x_n) = (\phi_{a_1}(x_1), \ldots, \phi_{a_n}(x_n))$, where $\phi_{a_i}$ is the partial conjugation outer automorphism associated above to the element $a_i$. Let $\{s_1, \ldots, s_n\}$ and $\{t_1, \ldots, t_n\}$ be subsets of $A$ with $n$ elements each. Let $G_{(s_1, \ldots, s_n)} = (S_1 \times \ldots \times S_n) \rtimes_{\phi_{(s_1, \ldots, s_n)}} \BZ$ and $G_{(t_i, \ldots, t_n)} = (S_1 \times \ldots \times S_n) \rtimes_{\phi_{(t_1, \ldots, t_n)}} \BZ$ be two semidirect products with such outer actions. Then $G_{(s_1, \ldots, s_n)}$ is isomorphic to $G_{(t_i, \ldots, t_n)}$ if and only if for the underlying sets $\{s_1, \ldots, s_n\} = \{t_1, \ldots, t_n\}$. 
\end{Lemma}

\begin{Proof}
($\Rightarrow$) Suppose there is an isomorphism $\theta$ between the two extensions.  Then $\theta$ must preserve the commutator subgroup, a characteristic group, so it induces an automorphism of $S_1 \times \ldots \times S_n$, say $\psi$. By Lemma \ref{lemstraightening-up-lemma}, $\theta$ must send each of the $n$ factors of $S_1 \times \ldots \times S_n$ in the domain isomorphically onto exactly one of the $n$ factors of $S_1 \times \ldots \times S_n$ in the range. Let $\sigma$ be the permutation from Lemma \ref{lemstraightening-up-lemma}.

Also, the associated map to $\theta$ on quotient groups must send the infinite cyclic quotient $G_{(s_1, \ldots, s_n)}/(S_1 \times \ldots \times S_n)$ isomorphically onto the infinite cyclic quotient $G_{(t_i, \ldots, t_n)}/(S_1 \times \ldots \times S_n)$. So, $\theta$ takes the generator, $z$, of $\BZ$ in $G_{(s_1, \ldots, s_n)}$ to the an element $cw^e$ in $G_{(t_i, \ldots, t_n)}$, where $c$ is some element of $S_1 \times \ldots \times S_n$ and $e$ is +1 or -1.

But also we know that $z$ centralises the factor $K_1$ in $S_i$, so its $\theta$-image $cw^e$ must centralize $\theta(K_1) = \psi(K_1) = Q$, say, in $S_{\sigma(i)}$, the copy of $S$ to which $S_i$ is sent under the isomorphism given by Lemma \ref{lemstraightening-up-lemma}, and act as conjugation by $t_{\sigma(i)}$ on $\theta(K_2) = \psi(K_2) = R$, say, in $S_{\sigma(i)}$.

Now, this implies that conjugation by $c_{\sigma(i)}$ has the same effect as conjugation by $w^{-e}$ on the subgroup $Q$ of $S_{\sigma(i)}$, which is isomorphic to $K_1$. By the Kurosh Subgroup Theorem, $Q = F*B_{1}^{d_1}*B_{2}^{d_2}* \ldots *B_{j}^{d_j}$, where $F$ is free, each $d \in S_{\sigma(i)}$ and $B_{j} \le K_{i_j}$, $i_j \in \{1,2\}$.

Since the Abelianization of $Q$ is trivial, we have that $F$ is trivial. As $Q \cong K_1$ and $K_1$ is freely indecomposable, we must have that $j = 1$. As $K_1$ is co-Hopfian, we must have that $B_{1} = K_{i_1}$.

So, $Q = \theta(K_1) = K_{i_1}^d$, with $d \in S_{\sigma(i)}$ and $i_1 \in \{1, 2\}$.

Suppose $\displaystyle \theta(K_1) = K_{1}^d$. Then

$$\displaystyle \theta(K_1^z) = \theta (K_1^z)$$
$$\displaystyle \theta(K_1) = \theta (K_1)^{\theta(z)}$$
$$\displaystyle K_{1}^d = (K_{1}^d)^{c_{\sigma(i)}w^e}$$
$$\displaystyle K_{1}^d = (K_{1}^{c_{\sigma(i)}w^e})^{d^{c_{\sigma(i)}w^e}}$$
$$\displaystyle K_{1}^d = (K_1^{c_{\sigma(i)}})^{d^{c_{\sigma(i)}t_{\sigma(i)}^e}}$$

So, $c_{\sigma(i)} = e$ and $d \in K_1$.

Suppose $\displaystyle \theta(K_1) = K_{2}^d$. Then

$$\displaystyle \theta(K_1^z) = \theta(K_1^z)$$
$$\displaystyle \theta(K_1) = \theta(K_1)^{\theta(z)}$$
$$\displaystyle K_{2}^d = (K_{2}^d)^{c_{\sigma(i)}w^e}$$
$$\displaystyle K_{2}^d = (K_{2}^{c_{\sigma(i)}w^e})^{d^{c_{\sigma(i)}w^e}}$$
$$\displaystyle K_{2}^d = (K_{2}^{c_{\sigma(i)}t_{\sigma(i)^e}})^{d^{c_{\sigma(i)}t_{\sigma(i)}^e}}$$

So, $\displaystyle c_{\sigma(i)} = t_{\sigma(i)}^{-e}$.

%Since $w^{-e}$ centralizes the $K_1$ factor in $S_{\sigma(i)}$ and acts as conjugation by $(t_{\sigma(i)})^{-e}$ on the $K_2$ factor, on the one hand, we have that $Q^{w^{-e}} = [(K_{i_1})^d]^{w^{-e}}$ is either $[(K_1)^d]^{w^{-e}} = (K_1^{w^{-e}})^(d^{w^{-e}}) = (K_1)^d = Q$ or $((K_2)^d)^{w^{-e}} = ((K_2)^d)^{t_{\sigma(i)}^{-e}} = [(K_2)^{t_{\sigma(i)}}]^{d^{t_{\sigma(i)}^{-e}}} = $.
%
%On the other hand, we have that $Q^{c_{\sigma(i)}} = = [(K_{i_1})^d]^{c_{\sigma(i)}}$ is either $(K_{1})^c_{\sigma(i)}^{d^{c_{\sigma(i)}}}$ or $((K_2)^d)^{c_{\sigma(i)}} = (K_{2})^c_{\sigma(i)}^{d^{c_{\sigma(i)}}}$
%
%In either event, we must have that $d \in K_1$.

By a symmetric argument, conjugation by $c_{\sigma(i)}$ has the same effect as conjugation by $w^{-e}$ after conjugation by $t_{\sigma(i)^e}$ on the subgroup $R$ of $S_{\sigma(i)}$, which is isomorphic to $K_2$. By the Kurosh Subgroup Theorem, $R = E*C_{1}^{f_1}*C_{2}^{f_2}* \ldots *C_{j}^{f_j}$, where $E$ is free, each $f \in S_{\sigma(i)}$ and $C_{k} \le K_{j_k}$, $j_k \in \{1,2\}$.

Since the Abelianization of $R$ is trivial, we have that $E$ is trivial. Also, as $R \cong K_2$ and $K_2$ is freely indecomposable and co-Hopfian, we must have that $R = \theta(K_2) = K_{j_1}^f$, with $f \in S_{\sigma(i)}$ and $j_1 \in \{1, 2\}$. Since $\theta$ is onto, the restriction of $\theta$ to $S_i$ must be onto. If $j_1 = i_1$, then the projection onto the $K_l$ factor of $\theta\bigr\rvert_{S_i}$, where $l \in \{1,2\} \backslash \{i_1\}$, would be trivial, and $\theta\bigr\rvert_{S_i}$ would not be onto. Thus, $j_1 \in \{1,2\} \backslash \{i_1\}$.

So, $R = \theta(K_2) = K_{j_1}^f$, with $d \in S_{\sigma(i)}$ and $j_1 \in \{1,2\} \backslash \{i_1\}$.

Suppose $\displaystyle \theta(K_2) = K_{1}^f$. Then

$$\displaystyle \theta(K_2^z) = \theta (K_2^z)$$
$$\displaystyle \theta(K_2^{s_i}) = \theta (K_2)^{\theta(z)}$$
$$\displaystyle \theta(K_{2})^{\theta(s_i)} = (K_{1}^f)^{c_{\sigma(i)}w^e}$$
$$\displaystyle (K_{1}^f)^{t_{\sigma(i)}^e} = (K_{1}^{c_{\sigma(i)}w^e})^{f^{c_{\sigma(i)}t_{\sigma(i)}^e}}$$
$$\displaystyle (K_{1}^{t_{\sigma(i)}^e})^{f^{t_{\sigma(i)}^e}} = (K_1^{c_{\sigma(i)}})^{f^{c_{\sigma(i)}t_{\sigma(i)}^e}}$$
$$\displaystyle (K_{1})^{f^{t_{\sigma(i)}^e}} = (K_1^{c_{\sigma(i)}})^{f^{c_{\sigma(i)}t_{\sigma(i)}^e}}$$

So, $c_{\sigma(i)} = e$.

Suppose $\displaystyle \theta(K_2) = K_{2}^f$. Then

$$\displaystyle \theta(K_2^z) = \theta(K_2^z)$$
$$\displaystyle \theta(K_{2}^{s_i}) = \theta(K_2)^{\theta(z)}$$
$$\displaystyle \theta(K_{2})^{\theta(s_i)} = (K_{2}^f)^{c_{\sigma(i)}w^e}$$
$$\displaystyle (K_{2}^f)^{t_{\sigma(i)}^e} = (K_{2}^{c_{\sigma(i)}w^e})^{f^{c_{\sigma(i)}w^e}}$$
$$\displaystyle (K_{2}^{t_{\sigma(i)}^e})^{f^{t_{\sigma(i)}^e}} = (K_{2}^{c_{\sigma(i)}t_{\sigma(i)}^e})^{f^{c_{\sigma(i)}t_{\sigma(i)}^e}}$$

So, $\displaystyle c_{\sigma(i)} = e$.

%So, $R = \theta(K_2) = K_{j_1}^f$, with $d \in S_{\sigma(i)}$ and $j_1 \in \{1, 2\} \backslash \{i_1\}$. Since $w^{-e}$ centralizes the $K_1$ factor in $S_{\sigma(i)}$ and acts as conjugation by $(t_{\sigma(i)})^{-e}$ on the $K_2$ factor, on the one hand, we have that $(R^{w^{-e}})^{t_{\sigma(i)^e}}$ is either $(((K_1)^f)^{t_{\sigma(i)}^{-e}}) = ((K_1^{t_{\sigma(i)}^{-e}})^{f^{t_{\sigma(i)}^{-e}}} = K_1^{f^{t_{\sigma(i)}^{-e}}}$ or $((K_2)^f)^{t_{\sigma(i)}} = (K_2^{t_{\sigma(i)}})^{f^{t_{\sigma(i)}^{-e}}}$.
%
%On the other hand, we have that, $R^{c_{\sigma(i)}} = = [(K_{j_1})^f]^{c_{\sigma(i)}}$ is either $[(K_{1})^f]^{c_{\sigma(i)}} = [(K_1)^{c_{\sigma(i)}}]^{f^{c_{\sigma(i)}}}$ or $((K_2)^f)^{c_{\sigma(i)}} = (K_{2})^c_{\sigma(i)}^{f^{c_{\sigma(i)}}}$
%
%In either event, w must have that $f \in K_1$.

But now, if $\theta(K_1) = K_1^d$ and $\theta(K_2) = K_2^f$, conjugation by $c_{\sigma(i)}$ must act trivially on $Q$ and trivially on $R$, which implies $c_{\sigma(i)}$ is trivial since $K_1$ and $K_2$, and hence $Q$ and $R$, are centerless; otherwise, if $\theta(K_1) = K_2^f$ and $\theta(K_2) = K_1^d$, conjugation by $c_{\sigma(i)}$ must act by conjugation by $t_{\sigma(i)}^{-e}$ on $Q$ and trivially on $R$.

Suppose for each $S_i$, $\theta(K_1) = K_2^f$ and $\theta(K_2) = K_1^d$. Let $x \ in K_2$ be an element other than the identity element and let $y = \theta(x) \in Q$. Then \\
$\theta(x^{z}) = \theta(x)^{\theta(z)}$ \\
$\theta(x^{s_i}) = y^{c_{\sigma(i)}t_{\sigma(i)}^e}$ \\
$\theta(x^{s_i}) = y$ as conjugation by $c_{\sigma(i)}$ must act by conjugation by $t_{\sigma(i)}^{-e}$ on $Q$ \\
$\theta(x^{s_i}) = \theta(x)$, which shows $\theta$ is not 1-1 for all $s_i \in A$ with $s_i$ not the identity, which contradicts the fact that $\theta$ is an isomorphism.

It follows that $\theta(K_1) = K_1^d$ and $\theta(K_2) = K_2^f$ and $c_{\sigma(i)}$ is trivial.

Thus $\theta$ takes $z$ to $w^e$. 

%Also $Q = K_1$, so $\psi$ preserves each $K_1$ factor.

Finally, if $s_i, t_{\sigma(i)} \in A$ and $x$ is any element of the $K_2$ factor in $S_i$, set $y = \theta(x) = \psi(x)$. Then

 $$ \theta(x^z) = \theta(x^{s_i}) = \theta(x)^{\theta(s_i)} = y^{\psi(s_i)} $$

while on the other hand,

  $$ \theta(x^z) = \theta(x)^{\theta(z)}
                      = y^w{^{e}}
                      = y^{(t_{\sigma(i)})^e} $$

and it follows that $y^{\psi(s_i)} = y^{(t_{\sigma(i)})^e}$, so that $\psi(s_i) = (t_{\sigma(i)})^e$; since no isomorphism of $K_1$ takes any given element of $A$ onto any other element of $A$ or the inverse of any other element of $A$, we must have $s_i = t_{\sigma(i)}^{\pm 1}$.

\indexspace

($\Leftarrow$) Suppose $\{s_1, \ldots, s_n\} = \{t_1, \ldots, t_n\}$. Choose $a \in G_{(s_1, \ldots, s_n)}$ with $aK_{(s_1, \ldots, s_n)}$ generating the infinite cyclic quotient $G_{(s_1, \ldots, s_n)}/K_{(s_1, \ldots, s_n)}$ and choose $b \in G_{(t_1, \ldots, t_n)}$ with $bK_{(t_1, \ldots, t_n)}$ generating the infinite cyclic quotient $G_{(t_1, \ldots, t_n)}/K_{(t_1, \ldots, t_n)}$. Set $\theta(a) = b$.

\indexspace

Send each element of $S_i$ (where $S_i$ uses an element $s_i$ in its semidirect product definition in the domain) to a corresponding generator of $S_i$ (where $S_i$ uses an element $t_i$ in its semidirect product definition in the range) under $\theta$.

\indexspace

Then $\theta: G_{(s_1, \ldots, s_n)} \rightarrow G_{(t_1, \ldots, t_n)}$ is a bijection. It remains to show $\theta$ respects the multiplication in each group.

\indexspace

Clearly, $\theta$ respects the multiplication in each $S_i$ and in $\BZ$

\indexspace

Finally, if $\alpha_i \in S_i$ and $a \in \BZ$,

\begin{tabular}{rcl}
$\theta(a \alpha_i)$            & = & $\theta(a)\theta(\alpha_i)$ \\
$\theta(\phi_{s_i}(\alpha_i)a)$ & = & $\phi_{t_i}(\theta(\alpha_i)) g(a)$
\end{tabular}

\indexspace

using the slide relators for each group and the fact that $s_i = t_i$, which implies $\phi_{s_i} = \phi_{t_i}$. So, $\theta$ respects the multiplication in each group. This completes the proof.

\end{Proof}

\begin{Lemma} \label{lemconderepimorphism-lemma}
Let $K$ and $A$ be as in the first paragraph of this section, let $(a_1, a_2, a_3, \ldots)$ be an enumeration of $A$, and let $(\omega, n) = (s_1, \ldots, s_n)$ and $(\eta, m) = (t_1, \ldots, t_m)$ be subsequences of $(a_1, a_2, a_3, \dots)$ (each with $n$ and $m$ distinct elements respectively for $n > m$).

\indexspace

Let $G_{(\omega, n)} = (S_1 \times \ldots \times S_n) \rtimes_{\phi_{(\omega, n)}} \BZ$ and $G_{(\eta, m)} = (S_1 \times \ldots \times S_m) \rtimes_{\phi_{(\eta, m)}} \BZ$ be two semidirect products. 	Then there is an epimorphism $g: G_{(\omega, n)} \rightarrow G_{(\eta, m)}$ if and only if $\{t_1, \ldots, t_m\} \subseteq \{s_1, \ldots, s_n\}$. 
\end{Lemma}

\begin{Proof}
The proof in this case is similar to the case $n = m$ in Lemma \ref{lemconderisomorphism-lemma}, except that the epimorphism $g$ must crush out $n-m$ factors of $K_{(\omega,n)} = S_1 \times \ldots \times S_n$ by Corollary \ref{corstraightening-up-corollary} and the Pidgeonhole Principle and then is an isomorphism on the remaining factors.

\indexspace

($\Rightarrow$) Suppose there is an epimorphism $g: G_{(\omega, n)} \rightarrow G_{(\eta, m)}$. Then $g$ must send the commutator subgroup of $G_{(\omega, n)}$ onto the commutator subgroup of $G_{(\eta, m)}$. By Corollary \ref{corstraightening-up-corollary}, $g$ must send $m$ factors of $K_{(\omega,n)} = S_1 \times \ldots \times S_n$ in the domain isomorphically onto the $m$ factors of $K_{(\eta,m)} = S_1 \times \ldots \times S_m$ in the range and sends the remaining $n - m$ factors of $K_{(\omega,n)}$ to the identity. Let $\{i_1, \ldots, i_m\}$ be the indices in $\{1, \ldots, n\}$ of factors in $K_{(\omega,n)}$ which are sent onto a factor in $K_{(\eta,m)} $ and let $\{j_1, \ldots, j_{n-m}\}$ be the indices in $\{1, \ldots, n\}$ of factors in $K_{(\omega,n)}$ which are sent to the identity in $K_{(\eta,m)} $. Then $g$ induces an isomorphism between $S_{i_1} \times \ldots \times S_{i_m}$ and $K_{(\eta,m)} $. Set $L_m = S_{i_1} \times \ldots \times S_{i_m}$

\indexspace

Also, by an argument similar to Lemma \ref{lemconderisomorphism-lemma}, $g$ sends sends the infinite cyclic group $G_{(\omega, n)}/K_{(\omega, n)}$ isomorphically onto the infinite cyclic quotient $G_{(\eta, m)}/K_{(\eta, m)}$.

\indexspace

Note that $L_m \rtimes_{\phi_{(s_{i_1}, \ldots, s_{i_m})}} \BZ$ is a quotient group of $G_{(\omega, n)}$ by a quotient map which sends $S_{j_1} \times \ldots \times S_{j_{n-m}}$ to the identity. Consider the induced map $g': L_m \rtimes_{\phi_{(s_{i_1}, \ldots, s_{i_m})}} \BZ \rightarrow G_{(\eta,m)}$. By the facts that $g'$ maps $L_m$ isomorphically onto $K_{(\eta,m)}$ and preserves the infinite cyclic quotients, we have that the kernel of $g$ must equal exactly $S_{j_1} \times \ldots \times S_{j_{n-m}}$; thus, by the First Isomorphism Theorem, we have that $g'$ is an isomorphism.

\indexspace

Finally, $g'$ is an isomorphism of $L_m \rtimes_{\phi_{(s_{i_1}, \ldots, s_{i_m})}} \BZ$ with $G_{(\omega, n)}$ which restricts to an isomorphism of $L_m$ with $S_{t_1} \times \ldots \times S_{t_m}$, so, as each $s_i$ and $t_i$ appears at most once, by an argument similar to Lemma \ref{lemconderisomorphism-lemma}, $\{t_1, \ldots, t_m\} \subseteq \{s_1, \ldots, s_n\}$.

\indexspace

($\Leftarrow$) Suppose $\{t_1, \ldots, t_m\} \subseteq \{s_1, \ldots, s_n\}$. Choose $a \in G_{(\omega, n)}$ with $aK_{(\omega, n)}$ generating the infinite cyclic quotient $G_{(\omega, n)}/K_{(\omega, n)}$ and choose $b \in G_{(\eta, m)}$ with $bK_{(\eta, m)}$ generating the infinite cyclic quotient $G_{(\eta, m)}/K_{(\eta, m)}$. Set $g(a) = b$.

\indexspace

Send each element of $S_i$ (where $S_i$ uses an element of order $t_i$ in its semidirect product definition in the domain) to a corresponding generator of $S_i$ (where $S_i$ uses an element of order $t_i$ in its semidirect product definition in the range) under $g$. Send the elements of all other $S_j$'s to the identity.

\indexspace

Then $g: G_{(\omega, n)} \rightarrow G_{(\eta, m)}$ is an epimorphism. Clearly, $g$ is onto by construction. It remains to show $g$ respects the multiplication in each group.

\indexspace

Clearly, $g$ respects the multiplication in each $S_i$ and in $\BZ$

\indexspace

Finally, if $\alpha_i \in S_i$ and $a \in \BZ$,

\begin{tabular}{rcl}
$g(a \alpha_i)$            & = & $g(a)g(\alpha_i)$ \\
$g(\phi_{s_i}(\alpha_i)a)$ & = & $\phi_{t_i}(g(\alpha_i)) g(a)$
\end{tabular}

\indexspace

using the slide relators for each group and the fact that $s_i = t_i$, which implies $\phi_{s_i} = \phi_{t_i}$. So, $g$ respects the multiplication in each group. This completes the proof.
\end{Proof}

Recall $K$ and $A$ are as in the first paragraph of this section, and $(a_1, a_2, a_3, \ldots)$ is an enumeration of $A$. Set $\Omega$ to be an uncountable set consisting of all subsequences of $(a_1, a_2, a_3, \ldots)$. For $\omega \in \Omega$ and $n \in \BN$, recall we have defined $(\omega, n)$ to be the finite sequence consisting of the first $n$ entries of $\omega$.

Set $G_{(\omega, n)} = (S \times S \times  \ldots \times S) \rtimes_{\phi_{(\omega, n)}} \mathbb{Z}$.

Proofs of Lemmas \ref{lemsessplitting-lemma} and \ref{lemladderdiagramsnonequiv-lemma} first appeared in \cite{Rolland1}; the reader is referred there for proofs.

\begin{Lemma} \label{lemsessplitting-lemma}
$G_{(\omega, n)} \cong  S \rtimes_{\phi_{a_{i_n}}} G_{(\omega, n-1)}$, where $\phi_{a_{i_n}}$ is partial conjugation by $a_{i_n}$.
\end{Lemma}

Now, this way of looking at $G_{(\omega, n)}$ as a semi-direct product of $S$ with $G_{(\omega, n-1)}$ yields an inverse sequence $(G_{(\omega, n)}, \alpha_n)$, which looks like

\begin{diagram}
G_{(\omega, 0)} & \lTo^{\alpha_{0}} & G_{(\omega, 1)} & \lTo^{\alpha_{1}} & G_{(\omega, 2)} & \lTo^{\alpha_{2}} & \ldots 
\end{diagram}

with bonding maps $\alpha_i: G_{(\omega, i+1)} \rightarrow G_{(\omega, i)}$ that each crush out the most recently added copy of $S$.

\indexspace

A subsequence will look like 

\begin{diagram}
G_{(\omega, n_0)} & \lTo^{\alpha_{n_0}} & G_{(\omega, n_1)} & \lTo^{\alpha_{n_1}} & G_{(\omega, n_2)} & \lTo^{\alpha_{n_2}} & \ldots 
\end{diagram}

with bonding maps $\alpha_{n_i}: G_{\omega, n_j)} \rightarrow G_{\omega, n_i)}$ that each crush out the most recently added $n_j - n_i$ copies of $S$.

\indexspace

\begin{Lemma} \label{lemladderdiagramsnonequiv-lemma}
If, for inverse sequences $(G_{(\omega, n)}, \alpha_{n})$ and $(G_{(\eta, m)}, \beta_{m})$, where $\alpha_{n}: G_{(\omega, n)} \rightarrow G_{(\omega, n-1)}$ is the bonding map crushing out the most recently-added copy of $S$ and $\beta_{m}: G_{(\eta, m)} \rightarrow G_{(\eta, m-1)}$ is the bonding map crushing out the most recently-added copy of $S$, $\omega$ does not equal $\eta$, then the two inverse sequences are not pro-isomorphic.
\end{Lemma}

\begin{RestateTheorem}{Theorem}{thmpseudo-collars}{Uncountably Many Pseudo-Collars on Closed Manifolds with the Same Boundary and Similar Pro-$\pi_1$}
Let $M^n$ be a closed smooth manifold ($n \ge 6$) with $\pi_1(M) \cong \BZ$ and let $S$ be the finitely presented group $K*K$ which is the free product of 2 copies of an admissible group (defined below) group $K$. Then there exists an uncountable collection of pseudo-collars $\{N^{n+1}_{\omega}\ |\ \omega \in \Omega\}$, no two of which are homeomorphic at infinity, with  $\partial N^{n+1}_{\omega} = M^n$, each with fundamental group at infinity that may be represented by an inverse sequence

\begin{diagram}[size=14.5pt]
\BZ & \lOnto^{\alpha_1} & G_{1} & \lOnto^{\alpha_2} & G_{2} & \lOnto^{\alpha_3} & G_{3} & \lOnto^{\alpha_4} & \ldots \\
\end{diagram}

with $\ker(\alpha_i) = S$ for all $i$.
\end{RestateTheorem}

\begin{Proof}
For each element $\omega \in \Omega$, the set of all increasing sequences of prime numbers, we will construct a pseudo-collar $N_{\omega}^{n+1}$ whose fundamental group at infinity is represented by the inverse sequence $(G_{(\omega, n)}, \alpha_{(\omega, n)})$. By Lemma \ref{lemladderdiagramsnonequiv-lemma}, no two of these pseudo-collars can be homeomorphic at infinity, and the Theorem will follow.

\indexspace

To form one of the pseudo-collars, start with a sequence $\omega \in \Omega$, set $s_i = \omega(i)$, let $M = \BS^1 \times \BS^{n-1}$ with fundamental group $\BZ$, and then blow $M$ up, using Theorem 1.1 of \cite{Rolland1}, to a cobordism $(W_{(s_1)}, M, M_{(s_1)})$ corresponding to the group $G_{(s_1)}$ ($s_1$ an element of $A$)..

\indexspace

We then blow this right-hand boundaries up, again using Theorem 1.1 of \cite{Rolland1} and Lemma \ref{lemsessplitting-lemma}, to cobordisms $(W_{(s_1, s_2)}, M_{(s_1)}, M_{(s_1, s_2)})$ corresponding to the group $G_{(s_1, s_2)}$ above.

\indexspace

We continue in the fashion \textit{ad infinitum}.

\indexspace

The structure of the collection of all pseudo-collars will be the set $\Omega$ described above.

\indexspace

We have shown that the pro-fundamental group systems at infinity of each pseudo-collar are non-pro-isomorphic in Lemma \ref{lemladderdiagramsnonequiv-lemma}, so that all the ends are non-diffeomorphic (indeed, non-homeomorphic).

\indexspace

This proves we have uncountably many pseudo-collars, each with boundary $M$, which have distinct ends.
\end{Proof}

\bibliography{02-more-examples-of-pseudocollars} %use this if you use bibtex 
\bibliographystyle{plain}

\end{document}